\input amstex
\input xy
\xyoption{all}
\documentstyle{amsppt}
\document
\magnification=1200
\NoBlackBoxes
\nologo
\hoffset1.5cm
\voffset2cm
\vsize15.5cm

\bigskip


\smallskip
\centerline{\bf  TIME AND PERIODICITY FROM PTOLEMY TO SCHR\"ODINGER: }

\medskip

\centerline{\bf  PARADIGM SHIFTS vs CONTINUITY}

\medskip

\centerline{\bf IN HISTORY OF MATHEMATICS}

\bigskip

\centerline{\bf Yuri I.~Manin}

\bigskip

\hfill{\it To Bob Penner, cordially}

\bigskip


{\bf Abstract.} I briefly consider the Kuhnian notion of ``paradigm shifts'' applied
to the history of mathematics and argue that the  succession and 
intergenerational continuity of  mathematical thought was undeservedly
neglected in the historical studies. To this end, I focus on the history
of mathematical theory of time and periodicity, from Ptolemy's epicycles to Schr\"odinger's 
quantum amplitudes interference and contemporary cosmological models.

\bigskip

\smallskip


\bigskip {\

\centerline{\bf Contents}

\smallskip

 Introduction
 
 1. Brief summary and plan of exposition

2. Mathematics and physics of periodicity

a. Antiquity: Euclidean geometry, Ptolemy's epicycles, Antikythera

$\quad\quad$ Digression 1: number $\pi$

b. Fourier sums and Fourier integrals: epicycles' calculus

$\quad\quad$Digression 2: number $e$ and ``computational consciousness''

c. Quantum amplitudes and their interference

 References

\newpage

\centerline{\bf  Introduction}

\medskip

In his influential treatise  [Ku70],  Thomas Kuhn  developed an approach
to the history of natural science(s) based upon the  assumption  
that this history can be naturally subdivided into periods.  According to Kuhn,
the  transitions from one period to the next one (called ``revolutions'') are characterised by 
a radical change of the basic assumptions, experimental and observational practices,
and acceptable types of argumentation. Any such set of assumptions
is shared by the learned community during each  development phase of 
``normal science'', and its change is called a ``paradigm shift''.

\medskip

Kuhn himself was reluctant about extending this view to the history,
philosophy, religion, and much of the social science(s). He believed that they are formed rather
  by a ``tradition of claims, counterclaims, and debates over fundamentals.''
 \medskip
 
 The motivation of  this  brief essay
 was a desire to discuss the applicability of the Kuhnian  view on history of
 mathematics. I  argue that  the  succession and 
intergenerational continuity of  mathematical thought was undeservedly
neglected in the sciencessical studies. To this end, I focus on the history
of the the mathematical theory of time and periodicity, from Ptolemy's epicycles to Schr\"odinger's 
quantum amplitudes interference and Feynman integrals.

\smallskip

According to the concise description in [Da09], my essay
lands somewhere in the uncharted territory between History of Science and Science Studies.
Kuhn's book originated Science Studies ``as a self--conscious field of inquiry'' ([Da09], p.~801.)
Hence this article belongs to it. But it focuses on the intrinsic continuity and the
peculiarities of forms of historical legacy in understanding 
space, time, and periodicity that, for many historians, might be completely
outside their fields of vision.

 \medskip
 
 If one rejects, as I do here,  the assumption about (this particular flow of) history as a sequence of revolutions, 
 then the idea of paradigm shifts cannot claim anymore its leading role.
 
 \medskip
 
 I accept here the more general viewpoint that Mathematics has a position mediating, or bridging, daily life, common sense,
 philosophy, and physics. Those fragments of mathematical knowledge that can become
 subjected to
 ``reality tests'' part are more sensitive to respective
 ``revolutions'' or ``paradigm shifts'', whereas those parts
 that are closer to ``pure mathematics'' show rather a kind of continuous
 development as is argued in this paper. Cf. also [Pa90], [Si12], [W10].
 
 \medskip
 
 Also, some light on my position can be thrown by comparison of the history of developing
 knowledge on the scale of civilisations with the history of development
of cognition in the individual brain of a growing human being
(cf. [MM17].)

\medskip

{\it Acknowledgements.} A project of this paper was conceived during
a brief stimulating  conversation of Yu.M. with Lorraine Daston about the role
of Kuhn's ``paradigm shifts'' doctrine in the history of applied (or rather, ``applicable'') mathematics
in a broad sense of this word. 

\smallskip
Andreea S.~Calude provided 
informative data about the prehistory of the  cognitive behaviour
of early  {\it Homo Sapiens.} 

\smallskip

When this paper was already written, Matilde Marcolli
drew my attention to the article [Ga01] which contains a 
very careful and sensible survey of a considerable part of the same
historical background (but excluding  quantum mechanics, cosmology and my remarks 
on ``computational consciousness'').

\smallskip

I am very grateful to them for inspiring communication.

\bigskip

\hfill{\it oritur sol et occidit et ad locum suum revertitur}

\hfill{\it lustrans universa in circuitu pergit spiritus et in circulos suos revertitur}

\smallskip

\hfill {\it VULGATA CLEMENTINA, Ecclesiastes 1:5--1:6}

\medskip

\centerline{\bf  1. Brief summary and plan of exposition}

\medskip

I will start with a few words about notions and formulas
summarising some basic {\it mathematical}  tools used in the {\it contemporary} discussions of time and periodicity.

\smallskip

Fundamental is the fact that these tools are subdivided into two complementary parts:
geometric ones involving {\it space/time intuition} (as in Euclid's {\it Elements})
and algebraic/calculus ones involving formulae and computations and generally having
linguistic character. Arguably (cf. [Ma15] and references therein) this is one reflection
of the general dynamic patterns of interactions between right and left brain.

\smallskip

Start with an Euclidean plane $\bold{P}$ endowed with Euclidean metric. 
Then a choice of a point $0$, of a line $L$ passing through it, and of its orientation,  determines
an identification of the set of all points of this line with
the set of real numbers $\bold{R}$: this is its ``coordinatisation''. Call it
the $x$--line and denote now $L_x$.

\smallskip

Choose now another line passing through the same point, oriented, and 
orthogonal to the $x$--line; call it $y$--line $L_y$.
\smallskip
Now we can construct a ``coordinatisation'' of the whole plane $\bold{P}$
i.e., the identification of the set of points of $\bold{P}$ with the set
of ordered pairs of real numbers $\bold{R}^2.$ 

\smallskip

At this point, we can start describing various figures, actors of plane Euclidean geometry,
 by equations and inequalities between various
algebraic expressions involving $x$ and $y$. So, for example, a circle of
radius $r$ whose center is a point $(x_0,y_0)$, is the set of all points
$(x,y)$ whose coordinates satisfy the equation  $(x-x_0)^2+ (y-y_0)^2=r^2.$

We may call it {\it Cartesian picture of geometry.}

\smallskip

Similarly, using 3 coordinates one gets an algebraic picture of Euclidean geometry of
space; passing to 4 coordinates, with time axis added to three space axes, we get the
scene for Newtonian mechanics. But some fragments of this scene were 
already recognisable in the world pictures going back to the times of Archimedes
and Ptolemy, as the celebrated Antikythera mechanism modelling the movement of heavens and
relating them to the chronological dating of historical events. ([J17]).

\smallskip

Arguably, one important  contribution of history of ``periodicity'' to mathematics was
the crystallisation of the notions of ``definition'', initially emerging as
secondary to the notions of ``axiom'' and ``theorem'' as in
Euclid's ``Elements''. 


\smallskip

The less obvious one  was a ``reification'' of the idea of symmetry: statements and proofs
of most theorems of Euclidean geometry are not dependent on the choice of
origin of coordinates and therefore invariant with respect to parallel shifts
of the whole space, and also with respect to rotations, conserving
angles. Thanks to this, one can introduce ``Cartesian'' coordinates also
on the space of all Euclidean symmetries of an Euclidean plane/space.

\smallskip

Finally, I must mention that, using the terminology of one of the schools
of Science Studies, when I briefly quote and/or interpret 
mathematical intuition and historical data, I appeal mostly to
{\it ``ethnomathematics in the European context'' } leaving aside many
interesting achievements and inputs that came from Eastern, Chinese  
and other regions of the global world. For a much more complete and balanced treatment, 
see e.~g. [Wo16].

\newpage

\centerline{\bf 2. Mathematics and physics of periodicity}

\bigskip

{\bf a. Antiquity: Euclidean geometry, Ptolemy's epicycles, Antikythera
Mechanism.}

\hfill{\it This book I bought in Venice for one ducat in the year 1507}
\smallskip
\hfill{Albrecht D\"urer inscription in Euclid's book}

\hfill{ from his library}

\medskip

Before starting the central themes of our discussion, I must say
explicitly that accumulation and intergenerational transmission
of knowledge, became possible only at a certain stage of development
of human language(s), and somewhat later, of  written languages.

\smallskip

Moreover, as I argued in [Ma07], pp.~159--167 and 169--189,
the most important new functions of emerging language consisted
{\it not} in the transmission of concrete information about
``here and now''  (``in this grove a deer is grazing''), but rather
in creation of {\it ``spaces of possibilities''}. Gods, heavens and netherworlds
powerfully influenced human's collective behaviour,  even if they could never
be located  here and now.

\smallskip

Since the concept of {\it here and now} itself later entered
physical theories as coordinate origin, it would be interesting
to trace its history as far back in time as possible. I am grateful
to Andreea Calude who informed me that deep reconstruction 
(to about  $15\cdot10^3$ years back from now)
seemingly recovers old common Indoeuropean  roots for ``now" but not
for ``here'', cf. [PagAtCaMe13]. Perhaps, a psychologically
motivated substitution for ``here'' was furnished by very old
(``ultraconserved'') words for ``I'' and ``you''.

\smallskip

Passing now to the real origins of modern scientific knowledge about the
Solar System and the Universe
in the Greco--Roman and Hellenistic worlds,  we see
that its foundations were laid between 300 BCE and 200 CE
and connected in particular with the names of Euclid, Archimedes, and Ptolemy.
The history of ``here and now'', however, must alert us to
the tracing also of the background history of the development
of various new ``languages of science'', of translations
and mutual interactions between these languages, and
their intergenerational functioning.
\smallskip

Euclid of Alexandria conjecturally lived and worked at about
325 BCE -- 265 BCE  in the south
Mediterranean  Greek colonial city. He created the richest
and at that time logically perfect axiomatic description of two-- and
three--dimensional spaces with metric and their symmetry groups that
were made explicit only many centuries  later when the language of
coordinates was created and one could speak about geometry
using languages of algebra/calculus etc. Still, perception of Euclid's
``Elements'' as {\it the} foundational, almost {\it sacral} treatise
survived till the 19th century: in particular new editions and translations 
of his {\it Elements} after spreading of printing were second only to the
Bible.

\smallskip
See a very remarkable book [By1847] by
Oliver Byrne, {\it ``surveyor of her Majesty's settlements in the Falkland Islands 
and author of numerous mathematical works''}, where he keeps texts
of all his geometric chapters but
rewrites all of Euclid's definitions (axioms), statements and proofs
in pseudo--algebraic formulas in which traditional for us
letters {\it a,d,c, $\dots$, x,y,z} serving as notations  (for us) variables, constants, functions etc.
are replaced by coloured pictures of angles, triangles et al.

\smallskip

Claudius Ptolemy conjecturally was born about 85 CE in Egypt and died
about 165 CE in Alexandria, Egypt. His greatest achievement  described in the
``Almagest'' is a dynamical model of the Solar System. This model is
geocentric. This is justified by the fact that all our observation of planets and Sun
are made from the Earth. It represents the visible movements 
of the planets and the sun as complex combination of uniform circular motions
along {\it epicycles}, whose centres also move uniformly along
their ``secondary'' epicycles, and finally various centres themselves 
are cleverly displaced from their expected ideal positions.
 
\smallskip

We do not know much about the computational devices that were used in antiquity
in order to make Ptolemy's model and other models of observable periodicities
such as lunar phases quantitatively comparable with
observations. However, one remarkable archaeological discovery
was made in 1900 when a group of sponger fishers from Greece during of bad weather
anchored their boats near the island Antikythera and while they were diving
discovered at a depth of 42 meters an ancient shipwreck. Besides 
bronze and marble statues, it contained a very corroded lump of bronze.
All these remnants were transferred to the National Archaeological Museum
in Athens, and after several decades of sophisticated studies and reconstructions,
a general consensus arose summarised in [Sp08] as follows:

\smallskip

{\it The Antikythera mechanism is an ancient astronomical calculator that contains
a lunisolar calendar, predicts eclipses, and indicates the moon's
position and phase. Its use of multiple dials and
interlocking gears eerily foreshadows modern computing concepts
from the fields of digital design, programming, and software engineering.}

\smallskip

For a description of continuing disagreements about details of the
reconstruction, see [Fr02] et al.
\medskip

{\bf Digression  1: the number $\pi$.} In the history of geometric models of
periodicity, the number $\pi$ plays a crucial role. Since Babylonian and Egyptian
times, $\pi$ was considered (``defined'') as the ratio of the length of a circle
to its diameter that can be measured in the same way as other physical
constants are measured. So in order to get an (approximate) value of $\pi$, one can
first say,  draw a circle using compasses, and then measure its length using
a string. Independence of the result on the diameter is also an
experimental fact which very naturally appears during land surveying.
Finally, the approximate values are always rational numbers, or rather,
{\it names} of some rational numbers, that can be transferred by
means acceptable in the relevant culture: see a very expressive
account by Ph.~E.~B.~Jourdain [Jou1956] written at the beginning of
the XXth century.

\smallskip
Arguably, the first modern approach to $\pi$ was found by 
Archimedes (about 287 -- 212 BCE). This approach consisted in
approximating $\pi$ from below by the values of perimeters
of inscribed regular $n$--gons (diameter is for simplicity taken as unit
of length). Manageable and fast converging formulas
for consecutive approximations are obtained by passing
from an $N$--gon to $2N$--gon etc.

\medskip

{\bf b. Fourier sums and Fourier integrals: epicycles' calculus.} As
we reminded in Sec.~1, after choosing orthogonal coordinates and scale
identifying an Euclidean plane with $\bold{R}^2$, one can describe
the circle of radius $r_0$ with centre $(x_0,y_0)$ as the set of
points $(x,y)$ such that  $(x-x_0)^2+(y-y_0)^2=r_0^2.$ 
The variables change $x=r_0(x_0+sin \,2\pi t)$, $y=r_0(y_0+ cos \,2\pi t)$
describes then the movement of a point along this circle,
with angular velocity one, if $t$ is interpreted as time flow. Replacing
$t$ by $v_0t$ we can choose another velocity.

\smallskip

In turn, we can put in the formulae above $x_0 = r_1(x_1+sin \,2\pi v_1t)$,
$y_0 = r_1(y_1+ cos \,2\pi v_1t)$, in order to make the centre $(x_0,y_0)$
move along another circle with uniform angular velocity, etc.
We get thus an analytic description of Ptolemy's picture, or rather its projection
on a coordinate plane in our space, which can be complemented by projections on
other planes. 
\smallskip

In order to use it for computational purposes, we must input the observable values
of $(x_i,y_i)$  and $v_i$, $i=0,1,\dots$, for, say, planetary movements.
The Antikythera mechanism served as a replacement of these formulae
for which the language was not yet invented and developed.
This language in its modern form 
and the analytic machinery were introduced only in the XVIII -- XIX centuries:
i.e.~Fourier sums/series $\sum_i (a_i sin\, i t + b_i cos\, i t)$
and more sophisticated Fourier integrals were initiated
by Jean--Baptiste Joseph Fourier (1768 -- 1830).

\smallskip

Joseph Fourier had a long and complicated social and political career starting with
education in the Convent of St.~Mark, and including service
in the local Revolutionary Committee during the French Revolution,
imprisonment during the Terror time,  travels with Napoleon to Egypt,
and  office of the Prefect of the Department of Is\`ere (where Joseph
Fourier was born).

\smallskip

Returning briefly to Fourier's mathematics, I would like to stress also
an analogy with Archimedes legacy, namely, observational astronomy
and mathematics of his ``Psammit'' (``The Sand Reckoner''). 
Archimedes wanted to estimate the size of the observable universe
giving an estimate of the number of grains of sand needed to fill it.
Among other difficulties he had to overcome, was the absence of
language (system of notation) for very large (in principle, as large as one wishes)
integers. He solved it by introducing inductively powers of 10,
so that any next power might be equal to the biggest number, defined at the previous step.
 
\medskip

{\bf Digression 2: the number $e$ and ``computational consciousness''.} The famous
Euler number $e=2,7182818284590\dots$ and his series
$$
e^x= 1+\sum_{n=1}^{\infty} \frac{x^n}{n!}
$$
were only the last steps of a convoluted history, with decisive contributions
due to  John Napier (1550--1617, Scotland), Henry Briggs (1561--1630, England),
Abraham de Moivre (1667--1754, France), among others, and finally Leonhard Euler
(1707, Basel, Switzerland -- 1783 St Petersburg, Russia).

\smallskip 

As already with Archimedes, and later with the
Masters of the Antikythera mechanism, one of the great motivations of the studies in this
domain was the necessity to devise practical tools for computations
with big numbers and/or numbers whose decimal notation included  many digits
 before/after the decimal point/comma: this is what I call here
 ``computational consciousness''. This ancient urge morphed now
 into such ideas as ``Artificial intellect'' and general identification
 of the activity of neural nets with computations.
 
 \smallskip
 
 So, for example, Briggs logarithm tables allowing
 to efficiently replace (approximate) multiplications by additions
 consisted essentially in the tables of numbers $10^7\cdot (1-!0^{-7})^N,
 N=1,2,3, \dots 10^7$. The future Euler's number $e$
 was hidden here as a result of passing to the limit
$e^{-1}= \roman{lim} (1- N^{-1})^N, N\to \infty ,$ which Briggs never
made explicit. However, the way Napier approached logarithms
included approximate calculations of logarithms of the function $\roman{sin}$,
which made be considered as the premonition of the Euler formula
$e^{ix}= \roman{cos}\,x + i\,\roman{sin}\,x$ that later played the key role
in mathematical foundations of quantum mechanics.
 
 \bigskip
 
{\bf c. Quantum amplitudes and their interference.} As we have seen, the basic scientific 
meta--notions
of ``observations'' and ``mathematical models'' explaining and predicting
results of observations, go back to deep antiquity.  The total body of scientific
knowledge accumulated since then, 
was enriched during the XIX--th and XX--th centuries also by recognition that ``scientific laws'',
that is, the central parts of mathematical models
explaining more or less directly  the results of observations, are
{\it qualitatively different} at various space/time scales: see a comprehensive survey
[`t HoVa14], in particular, the expressive table on pp.~100--101.

\smallskip

A breakthrough in understanding physics at the very large scale Universe (cosmology)
was related to Einstein's general relativity (or gravity) theory, whereas
on the very small scale, the respective breakthrough came with quantum mechanics
and later quantum field theory. Bridging these two ways of understanding
Nature still remains one of the main challenges for modern science.

\smallskip
One can argue that  an ``observable'' bridge between these two
scales is the existence and cognitive activities of {\it Homo Sapiens}  
on our Earth (and possibly elsewhere), but
the discussion of the current stage of ``observations'' and ``explanations''
in biology would have taken us too far away from the subject of this short essay;
cf.~[MM17]. Anyway, the key idea of scientific observation includes
some understanding of how a subject of human scale can interact
with objects of cosmic/micro scales.

\smallskip

Studying the small scale physics unavoidably involved the necessity
of working out mathematical models of probabilistic behaviour
of elementary particles that was observed and justified 
in multiple experiments.
It was preceded by a remarkable cognitive passage: from the observable
properties of chemical reactions to Mendeleev's intellectual construct of the
Periodic Table to the images of atoms of the Chemical Elements 
as analogs of the Solar System with nucleus for Sun
and electrons orbiting like planets. This cognitive passage might be
compared with the evolution of astronomy from antiquity to Copernicus, Galileo and Newton.

\smallskip

When experimental methods were developed for working quantitatively
with unstable (radioactive) atoms, small groups of electrons,
etc., a new theoretical challenge emerged: observable data
involved random, probabilistic behaviour, but the already well 
developed mathematical tools for describing randomness
did not work correctly in the microworld!

\smallskip

The emergent quantum mechanics postulated that ``probabilities''
of classical statistics, expressed by real numbers between 0 and 1,
must be replaced in the microworld by {\it probability amplitudes},
whose values are {\it complex numbers} that after some
normalisation become complex numbers lying in the complex plane
on the circle of radius 1 and centre 0. It must be then explained
how to pass from the hidden quantum mechanical picture to the observable
classical statistical picture. Many different paths along this thorny way
were discovered in the 20's of the XX century, in particular,
in classical works of Werner Heisenberg, Wolfgang Pauli, Erwin
Schr\"odinger, et al. One way of looking at ``quantization'' of
the simplest classical system, point--like body in space,
is this. Any classical trajectory of this system is a curve
in the space of pairs {\it (position, momentum)}.
Here position and momentum can be considered as Cartesian coordinate
triples whose values, of course, pairwise commute.
To the contrary, in the quantum mechanical space, the commutator between
position and momentum is not zero. This can be envisioned
as a replacement of possible classical trajectories of such a system
by their wave functions which are not localised.
Probabilistic data occur when one adds, say, one more point--like body
and/or interaction with a macroscopic environment
created by an experimenter.

\smallskip

Mathematical descriptions of all this are multiple and all
represent a drastic break with ``high school'', or ``layman'',
intuition. One remarkable example of pedagogical difficulties
of quantum physics can be glimpsed in the famous
Lectures on Quantum Mechanics by the great Richard Feynman.

\smallskip

In our context, the most essential is the fact that quantum interaction
in the simplest cases of quantum mechanics is described via
Fourier sums, series, and integrals in (finite dimensional)
{\it complex} spaces endowed with {\it Hermitean} metrics, in place of 
real Euclidean space with real metric.

\smallskip

The quantum mechanical amplitudes are given by
Fourier sums or series  of the form $\sum_n a_n e^{it}$
where $a_n$ are complex numbers, and $t$ is time,
whereas probabilities in classical statistic descriptions
are given by the similar sums with real $a_n$, and $it$
replaced by inverse (also real) temperature $-1/T$.

\smallskip

In this sense, quantum
mechanics
is a complexification of Ptolemy's epicycles.

\smallskip

In the currently acceptable picture,
our evolving Universe can be dissected into ``space sections''
corresponding to the values of global {\it cosmological}
time (e.g. in the so called Bianchi cosmological models)
to each of which a specific temperature of background
cosmic radiation can be ascribed. Going back in time, our
Universe becomes hotter, so that at the moment of the Big Bang (time = $0$)
its temperature becomes infinite. This provides
a highly romantic interpretation of the
correspondence  $-1/T \leftrightarrow it$.

\bigskip


\centerline{\bf References}

\medskip

[By1847] O.~Byrne. {\it The first six books of the Elements of Euclid in which coloured
diagrams and symbols are used instead of letters for the greater ease
of learners, By Oliver Byrne surveyor of her Majesty's settlements 
in the Falkland Islands and author of numerous mathematical works. London William
Pickering 1847}. 

Facsimile publication by Taschen.

\smallskip

[Da09] L.~Daston. {\it Science Studies and the History of Science.} 
Critical Inquiry, Vol.~35, No.~4, pp.~798--813. U.~of Chicago Press Journals, Summer of 2009.
https://doi.org/10.1086/599584

\smallskip

[Fr02] T.~Freeth. {\it The Antikythera mechanism.} Mediterranean Archaeology and
Archaeometry, vol.~2, no.~1, pp.~21--35.

\smallskip

[Ga01] G.~Gallavotti. {\it Quasi periodic motions from Hipparchus to Kolmogorov.}
Rendiconti Acc. dei Lincei, Matematica e Applicazioni., vol.~12, 2001, pp.~125--152.

\smallskip
[Gr17] St.~Greenblatt. {\it The rise and fall of Adam and Eve.} W.~W.~Norton \& Company, 2017.

\smallskip

['t HoVa14] G.~'t Hooft, St.~Vandoren. {\it Time in powers of ten. Natural Phenomena
and Their Timescales.} World Scientific, 2014.

\smallskip

[J17] A.~Jones. {\it A portable cosmos: revealing the Antikythera mechanism, scientific wonder
of the ancient world.} Oxford UP, 2017.

\smallskip

[Jou1956] Ph.~E.~B.~Jourdain. {\it The Nature of Mathematics.} Reproduced in the anthology
``The World of Mathematics'' by James R. Newman, vol. 1, Simon and Schuster, NY 1956.

\smallskip
[Ku70] Th.~Kuhn. {\it The Structure of Scientific Revolutions (2nd, enlarged ed.)}. 
University of Chicago Press, 1970. ISBN 0-226-45804-0.

\smallskip

[Ma15] Yu.~I.~Manin.  {\it De Novo Artistic Activity, Origins of Logograms, and Mathematical
Intuition. } In: Art in the Life of Mathematicians, Ed. Anna Kepes Szemer\'edi, AMS, 2015, pp. 187--208.

\smallskip

[Ma07] Yu. Manin. {\it Mathematics as Metaphor. Selected Essays, with Foreword
by Freeman Dyson.}

\smallskip
[MM17]  D.~Yu.~Manin, Yu.~I.~Manin. {\it Cognitive networks: brains, internet, and civilizations.}.
 In: Humanizing Mathematics and its Philosophy,  ed. by Bh.~Sriraman, Springer International
Publishing AG, 2017,  pp.~85--96. DOI 10.1007/978-3-319-61231-$7_{-}$9. arXiv:1709.03114

\smallskip

[PagAtCaMe13] M.~Pagel, Q.~D.~Atkinson, A.~S.~Calude, A.~Meade.
{\it Ultraconserved words point to deep language ancestry across Eurasia.}

www.pnas.org/cgi/doi/10.1073/pnas.1218726110

\smallskip

[Pa90]  D.~Park. {\it The How and the Why. An essay on the origins and development
of physical theory.} Princeton UP, New Jersey, 1990.

\smallskip


\smallskip

[Si12] K.~Simonyi. {\it A Cultural History of Physics. Translated by D.~Kramer.}
CRC Press, Boca Raton, London, New York, 2012.

\smallskip

[Sp08] D.~Spinellis. {\it The Antikythera Mechanism: a Computer Science Perspective.}
Computer, May 2008, pp.~22--27.

\smallskip

[W10] G.~Wolf.  {\it Orrery and Claw.} London Review of Books, 18 Nov. 2010, pp. 34--35.

\smallskip

[Wo16] D.~Wootton. {\it The Invention of Science: A New Histrory of the Scientific Revolution.}
Penguin, 1976, 784 pp.

\bigskip

{\bf Max--Planck--Institute for Mathematics, 

Vivatsgasse 7, Bonn 53111, Germany.}

\smallskip

manin\@mpim-bonn.mpg.de

\medskip

\enddocument